\documentclass{article}
\usepackage{graphicx}
\usepackage{psfrag}
\usepackage{amssymb}
\usepackage{amsmath}
\usepackage{amsthm}
\usepackage{amscd}

\newtheorem{thm}{Theorem}[section]
\newtheorem{prop}[thm]{Proposition}

\newtheorem{lem}[thm]{Lemma}

\newtheorem*{rh}{Riemann Hypothesis}

\newtheorem*{the condition}{The Sufficient Condition for the Riemann Hypothesis}

\theoremstyle{definition}
\newtheorem{defn}[thm]{Definition}

\title{Ordinality and Riemann Hypothesis I}
\author{Young Deuk Kim
\\SNU College\\
Seoul National University\\
Seoul 08826, Korea
\\(yk274@snu.ac.kr)}
\date{\today}

\begin{document}
\maketitle

\begin{abstract}
We present a sufficient condition for the Riemann hypothesis.
This condition is the existence of a special ordering on the set of finite products of distinct odd primes.

\vspace{.3cm}
\end{abstract}

\section{Introduction}

The zeta function
$$\zeta(s)=\sum_{k=1}^\infty\frac{1}{k^s}$$
was introduced by Euler in 1737 for real variable $s>1$.
In 1859, Riemann \cite{riemann} extended the function to the complex meromorphic function $\zeta(z)$
with only a simple pole at $z=1$, and
$$\zeta(z)=\sum_{k=1}^\infty\frac{1}{k^z}$$
on ${Re\,}z>1$.
\begin{thm}[\cite{SS}]\label{thm;zeta}
The zeta function has a meromorphic continuation into the entire complex plane, whose only
singularity is a simple pole at $z=1$.
\end{thm}

The zeta function has infinitely many zeros, but there is no zero in the region $\operatorname{Re}(z)\geq 1$.
\begin{thm}[\cite{RJ},\cite{SS}]
The only zeros of the $zeta$ function outside the critical strip
$0<\operatorname{Re}(z)<1$ are at the negative even integers, $-2,\ -4,\ -6,\ \cdots$.
\end{thm}

The most famous conjecture on the zeta function is the Riemann hypothesis.
\begin{rh}[\cite{bombieri},\cite{sarnak}]
The zeros of $\zeta(z)$ in the critical strip lie on the critical line $\operatorname{Re}(z)=\frac{1}{2}$.
\end{rh}

Suppose that $x$ and $y$ are real numbers with $0<x<1$. It is known that if $x+yi$ is a zero of the zeta function, then so are $x-yi$, $(1-x)+yi$ and $(1-x)-yi$.\\

\indent
Riemann himself showed that if $0\leq y\leq 25.02$ and $x+yi$ is a zero of the zeta function, then $x=\frac{1}{2}$.
Therefore the Riemann hypothesis is true up to height 25.02.
In 1986, van de Lune, te Riele and Winter \cite{LRW} showed that the Riemann hypothesis is true up to height 545,439,823,215. Furthermore, in 2021, Dave Platt and Tim Trudgian \cite{PT} proved that the Riemann hypothesis is true up to height $3\cdot 10^{12}$.\\

\indent
Therefore, to prove the Riemann hypothesis, it is enough to show that if $\frac{1}{2}<x<1$ and $y>0$, then $x+yi$ is not a zero of the zeta function.
In this paper, we study a sufficient condition for the Riemann hypothesis.
This condition is the existence of a special ordering on the set of finite products of distinct odd primes.
This condition inspired the author to propose a complete proof \cite{kim1,kim2} of the Riemann hypothesis.

\section{Preliminary Lemmas and Theorems}

The eta function
$$\eta(z)=\sum_{k=1}^\infty\frac{(-1)^{k-1}}{k^z}=1-\frac{1}{2^z}+\frac{1}{3^z}+\cdots$$
is convergent on $\operatorname{Re}(z)>0$, where we assume that $(-1)^0=1$ for the sake of simplicity.

\begin{thm}[\cite{broughan}]\label{thm:ext}
For $0<\operatorname{Re}(z)<1$, we have
$$\zeta(z)=\frac{1}{1-2^{1-z}}\eta(z).$$
\end{thm}

\noindent
The zeros of $1-2^{1-z}$ are on $\operatorname{Re}(z)=1$.
Therefore, in the critical strip $0<\operatorname{Re}(z)<1$, any zero of $\zeta(z)$ is a zero of $\eta(z)$.

\begin{lem}\label{lem1}
Let $0<x<1$ and $y\in\mathbb{R}$.
If $x+yi$ is a zero of $\zeta(z)$ then
\begin{equation*}
\sum_{k=1}^\infty \frac{(-1)^{k-1}}{k^x}\cos(y\ln k)=\sum_{k=1}^\infty \frac{(-1)^{k-1}}{k^x}\sin(y\ln k)=0.
\end{equation*}
\end{lem}

\begin{proof}
\begin{eqnarray*}
\frac{1}{k^{x+yi}}=k^{-x-yi}&=&e^{(-x-yi)\ln k}\\
&=&e^{-x\ln k}\left(\cos(y\ln k)-i\sin(y\ln k)\right)\\
&=&\frac{1}{k^x}\left(\cos(y\ln k)-i\sin(y\ln k)\right)
\end{eqnarray*}
Therefore, it follows directly from Theorem \ref{thm:ext}.
\end{proof}

\begin{lem}\label{lem2}
Let $0<x<1$ and $y\in\mathbb{R}$.
If $x+yi$ is a zero of $\zeta(z)$ then
\begin{equation*}
\sum_{k=1}^\infty \frac{(-1)^{k-1}}{k^x}\cos(y\ln (ak))=0
\end{equation*}
for all $a>0$ and
\begin{equation*}
\sum_{k=1}^\infty \frac{(-1)^{k-1}}{k^x}\sin(y\ln (bk))=0
\end{equation*}
for all $b>0$.
\end{lem}
\begin{proof}
\begin{eqnarray*}
\cos(y\ln (ak))&=&\cos(y\ln a+y\ln k)\\
&=&\cos(y\ln a)\cos(y\ln k)-\sin(y\ln a)\sin(y\ln k)
\end{eqnarray*}
\begin{eqnarray*}
\sin(y\ln (bk))&=&\sin(y\ln b+y\ln k)\\
&=&\sin(y\ln b)\cos(y\ln k)+\cos(y\ln b)\sin(y\ln k)
\end{eqnarray*}
Therefore, it follows directly from Lemma \ref{lem1}.
\end{proof}

\begin{lem}\label{lem3}
Let $0<x<1$ and $y\in\mathbb{R}$.
Suppose that $x+yi$ is a zero of $\zeta(z)$ and $q\geq 1$ is an odd number. Then
\begin{equation*}
\sum_{m=1}^\infty \frac{(-1)^{mq-1}}{(mq)^x}\cos(y\ln (mq))=0
\end{equation*}
and
\begin{equation*}
\sum_{m=1}^\infty \frac{(-1)^{mq-1}}{(mq)^x}\sin(y\ln (mq))=0.
\end{equation*}

\end{lem}
\begin{proof}
Since $q$ is an odd number, $(-1)^{mq-1}=(-1)^{m-1}$.
Therefore, from Lemma \ref{lem2}, we have
\begin{eqnarray*}
\sum_{m=1}^\infty \frac{(-1)^{mq-1}}{(mq)^x}\cos(y\ln (mq))&=&\sum_{m=1}^\infty \frac{(-1)^{m-1}}{(mq)^x}\cos(y\ln (mq))\\
&=&\frac{1}{q^x}\sum_{m=1}^\infty \frac{(-1)^{m-1}}{m^x}\cos(y\ln (mq))=0
\end{eqnarray*}
and
\begin{eqnarray*}
\sum_{m=1}^\infty \frac{(-1)^{mq-1}}{(mq)^x}\sin(y\ln (mq))&=&\sum_{m=1}^\infty \frac{(-1)^{m-1}}{(mq)^x}\sin(y\ln (mq))\\
&=&\frac{1}{q^x}\sum_{m=1}^\infty \frac{(-1)^{m-1}}{m^x}\sin(y\ln (mq))=0.
\end{eqnarray*}
\end{proof}

\begin{lem}\label{lemf}
If $0<x<1$ and $y\in\mathbb{R}$, then
\begin{equation*}
1-\sum_{k=1}^\infty \frac{1}{2^{kx}}\cos(ky\ln 2)+i\sum_{k=1}^\infty \frac{1}{2^{kx}}\sin(ky\ln 2)\neq 0
\end{equation*}
\end{lem}
\begin{proof}

Since $0<x<1$, we have
\begin{eqnarray*}
&& 1-\sum_{k=1}^\infty \frac{1}{2^{kx}}\cos(ky\ln 2)+i\sum_{k=1}^\infty \frac{1}{2^{kx}}\sin(ky\ln 2)\\
&=& 2-\sum_{k=0}^\infty \frac{\cos(ky\ln 2)-i\sin(ky\ln 2)}{2^{kx}}\\
&=& 2-\sum_{k=0}^\infty \frac{e^{-iky\ln 2}}{2^{kx}}\\
&=& 2-\sum_{k=0}^\infty \left(\frac{e^{-iy\ln 2}}{2^x}\right)^k\\
&=& 2-\frac{1}{1-\frac{e^{-iy\ln 2}}{2^x}}\\
&=&2-\frac{2^x}{2^x-e^{-iy\ln 2}}\\
&=& \frac{2^x-2e^{-iy\ln 2}}{2^x-e^{-iy\ln 2}}\\
&\neq& 0.
\end{eqnarray*}
\end{proof}

\noindent
Lemma \ref{lemf} can be restated as the following theorem.
\begin{thm}\label{thm:II}
For each $k\in\mathbb{N}$, let
$$\varphi(k)=\frac{(-1)^{k-1}}{k^{x+iy}},$$
where we assume that $(-1)^0=1$ for the sake of simplicity.
If $0<x<1$ and $y\in\mathbb{R}$, then we have
\begin{equation*}
\sum_{\ell=0}^\infty\varphi(2^\ell)\neq 0.
\end{equation*}
\end{thm}

\begin{proof}
Since $0<x<1$, we have
\begin{eqnarray*}
&&\sum_{\ell=0}^\infty\varphi(2^\ell)=1-\sum_{\ell=1}^\infty\frac{1}{(2^\ell)^{x+iy}}=1-\sum_{\ell=1}^\infty\frac{e^{-i\ell y\ln 2}}{2^{\ell x}}\\
&&\qquad\quad =1-\sum_{\ell=1}^\infty \left(\frac{e^{-iy\ln 2}}{2^x}\right)^\ell = 1-\frac{e^{-iy\ln 2}}{2^x-e^{-iy\ln 2}}=\frac{2^x-2e^{-iy\ln 2}}{2^x-e^{-iy\ln 2}}\neq 0.
\end{eqnarray*}
\end{proof}

\section{The Sufficient Condition for\\ the Riemann Hypothesis}

Let $\mathbb{N}$ be the set of natural numbers and $\mathbb{N}_0=\mathbb{N}\cup\{ 0\}$.

\begin{defn}
Let $Q$ be the set of finite products of distinct odd primes.
\begin{equation*}
Q=\{ p_1p_2 \cdots p_n \mid p_1, p_2,\cdots, p_n \mbox{ are distinct odd primes},\ n\in\mathbb{N}\}
\end{equation*}
For each $q= p_1p_2 \cdots p_n\in Q$, we define
$$\mbox{sgn\;}q=(-1)^n$$
where $p_1, p_2,\cdots, p_n$ are distinct odd primes.
\end{defn}

There are infinitely many orderings of $Q$.
\begin{defn}\label{def:qh}
Choose an ordering on $Q$ and let $$Q=\{q_1,q_2,q_3,q_4,q_5,\cdots \}.$$
\end{defn}

\begin{defn}
For $i,k\in\mathbb{N}$, let
\[\delta(k,i)=\left\{
\begin{array}{cl}
1 & \mbox{ if } k \mbox{ is a multiple of } q_i \\
0& \mbox{ otherwise}
\end{array}\right.
\]
and
$$f(k,h)=\sum_{i=1}^h(\mbox{sgn\;}q_i)\delta(k,i),\qquad f(k)=\lim_{h\to\infty}f(k,h)=\sum_{i=1}^\infty(\mbox{sgn\;}q_i)\delta(k,i).$$
\end{defn}

\noindent
Note that, for each $k$, there exist only finitely many $i$ such that $\delta(k,i)\neq 0$.

\begin{defn}
Suppose that $\frac{1}{2}<x<1$, $y>0$ and $x+yi$ is a zero of $\zeta(z)$.
For $k\in\mathbb{N}$, let
$$a_k=\frac{(-1)^{k-1}}{k^x}\cos(y\ln k), \qquad b_k=\frac{(-1)^{k-1}}{k^x}\sin(y\ln k),$$
where we assume that $(-1)^0=1$ for the sake of simplicity.
\end{defn}

By Lemma \ref{lem1}, we have
\begin{equation}\label{eq:zero1}
\sum_{k=1}^\infty a_k=\sum_{k=1}^\infty b_k=0.
\end{equation}
From Lemma \ref{lem3}, we have
\begin{equation*}
\sum_{m=1}^\infty a_{mq_i}=0\qquad \mbox{for all }q_i\in Q
\end{equation*}

\noindent
and therefore
\begin{equation}\label{eq:zero2}
\sum_{m=1}^\infty (\mbox{sgn\;}q_i)a_{mq_i}=0\qquad \mbox{for all }q_i\in Q.
\end{equation}

\begin{defn}
For $n,h\in\mathbb{N}$, let
$$C(n,h)=\sum_{k=1}^n \sum_{i=1}^h (\mbox{sgn\;}q_i)\delta(k,i)a_k=\sum_{k=1}^n f(k,h)a_k $$
and
$$S(n,h)=\sum_{k=1}^n \sum_{i=1}^h (\mbox{sgn\;}q_i)\delta(k,i)b_k=\sum_{k=1}^n f(k,h)b_k. $$
\end{defn}

\begin{prop}\label{prop:nh1}
For each $h\in\mathbb{N}$, we have
$$\lim_{n\to\infty} C(n,h)=\lim_{n\to\infty} S(n,h)=0$$
and therefore
$$\lim_{h\to\infty} \lim_{n\to\infty} C(n,h)=\lim_{h\to\infty} \lim_{n\to\infty} S(n,h)=0.$$
\end{prop}

\begin{proof}
From eq. (\ref{eq:zero2}), for all $q_i\in Q$, we have
$$\sum_{k=1}^\infty  (\mbox{sgn\;}q_i)\delta(k,i) a_k=\sum_{m=1}^\infty  (\mbox{sgn\;}q_i) a_{mq_i}=0.$$
Therefore
\begin{eqnarray*}
0&=&\sum_{i=1}^h\sum_{m=1}^\infty  (\mbox{sgn\;}q_i) a_{mq_i} \\
&=&\sum_{i=1}^h\sum_{k=1}^\infty  (\mbox{sgn\;}q_i)\delta(k,i) a_k \\
&=&\sum_{k=1}^\infty\sum_{i=1}^h  (\mbox{sgn\;}q_i)\delta(k,i) a_k \\
&=&\lim_{n\to\infty}\sum_{k=1}^n\sum_{i=1}^h  (\mbox{sgn\;}q_i)\delta(k,i) a_k\\
&=&\lim_{n\to\infty} C(n,h).
\end{eqnarray*}
In the same way, we have
$$\lim_{n\to\infty} S(n,h)=0.$$
\end{proof}

\begin{defn}
Let
$$\Gamma=\{2^\ell\mid \ell\in\mathbb{N}_0\}.$$
\end{defn}

\begin{lem}\label{lem:fk}
Recall
$$ f(k)=\sum_{i=1}^\infty(\mbox{sgn\;}q_i)\delta(k,i).$$
For all $k\in\mathbb{N}$, we have
\begin{eqnarray*}
f(k)&=&\left\{
\begin{array}{cl}
0& \mbox{if } k\in\Gamma\\
-1 & \mbox{otherwise}
\end{array}\right.
\end{eqnarray*}
\end{lem}

\begin{proof}
If $k\in\Gamma$, then $k$ is not a multiple of any element in $Q$. Therefore
$\delta(k,i)=0$ for all $i$ and hence $f(k)=0$.

\indent
Suppose that $k\notin\Gamma$ and
$$k=2^mp_1^{m_1}p_2^{m_2}\cdots p_n^{m_n},\quad m_1, m_2,\cdots, m_n\geq 1,\quad m\geq 0,\quad n\geq 1$$
is the prime factorization of $k$, where $p_1,p_2,\cdots, p_n$ are distinct odd prime divisors of $k$.
We have
\begin{eqnarray*}
&&\{q_i\in Q\mid \delta(k,i)=1\}\\
&=&\{p_1,\cdots,p_n,\ p_1p_2,\cdots,p_{n-1}p_n,\ p_1p_2p_3,\cdots,p_1p_2\cdots p_n\}.
\end{eqnarray*}
Therefore
$$f(k)=-\binom{n}{1}+\binom{n}{2}-\cdots +(-1)^n\binom{n}{n}=-1.$$
\end{proof}

Notice that, for each $k$,
 $$\sum_{i=1}^\infty (\mbox{sgn\;}q_i)\delta(k,i)a_k\quad\mbox{and}\quad \sum_{i=1}^\infty (\mbox{sgn\;}q_i)\delta(k,i)b_k$$
become finite sums because $\delta(k,i)=0$ except finitely many $i$. For each $n$, from Lemma \ref{lem:fk}, we have
\begin{eqnarray*}
\lim_{h\to\infty} C(n,h)&=&\lim_{h\to\infty}\sum_{k=1}^n \sum_{i=1}^h (\mbox{sgn\;}q_i)\delta(k,i)a_k\\
&=&\sum_{k=1}^n \sum_{i=1}^\infty (\mbox{sgn\;}q_i)\delta(k,i)a_k\\
&=&\sum_{k=1}^n f(k)a_k\\
&=&\sum^{1\leq k\leq n}_{k\notin\Gamma} (-a_k)
\end{eqnarray*}

\noindent
In the same way we have
$$\lim_{h\to\infty} S(n,h)=\sum^{1\leq k\leq n}_{k\notin\Gamma} (-b_k).$$

Therefore we have the following proposition.
\begin{prop}\label{prop:nh2}
$$\lim_{n\to\infty}\lim_{h\to\infty} C(n,h)=\sum_{k=1}^\infty f(k)a_k$$
$$\lim_{n\to\infty}\lim_{h\to\infty} S(n,h)=\sum_{k=1}^\infty f(k)b_k$$
\end{prop}

\indent
Up to now, we have worked with an arbitrary ordering on $Q$. To prove Riemann hypothesis, we need a special ordering on $Q$. If the following condition is true, we can prove the Riemann hypothesis.
\begin{the condition}
There exists an ordering on $Q$ such that
\begin{equation}\label{eq:condition1}
\lim_{n\to\infty}\lim_{h\to\infty} C(n,h)=\lim_{h\to\infty}\lim_{n\to\infty} C(n,h)
\end{equation}
and
\begin{equation}\label{eq:condition2}
\lim_{n\to\infty}\lim_{h\to\infty} S(n,h)=\lim_{h\to\infty}\lim_{n\to\infty} S(n,h).
\end{equation}
\end{the condition}

\begin{thm}
If the above condition is satisfied, then the Riemann hypothesis is true.
\end{thm}

\begin{proof}
Suppose that there exists an ordering on $Q$ satisfying eq. (\ref{eq:condition1}) and  (\ref{eq:condition2}).
Let $\frac{1}{2}<x<1$, $y>0$ and $x+yi$ is a zero of $\zeta(z)$. This leads to a contradiction.

\indent
From Proposition \ref{prop:nh1} and Proposition \ref{prop:nh2}, we have
$$\sum_{k=1}^\infty f(k)a_k=\lim_{n\to\infty}\lim_{h\to\infty} C(n,h)=\lim_{h\to\infty}\lim_{n\to\infty} C(n,h)=0$$
and
$$\sum_{k=1}^\infty f(k)b_k=\lim_{n\to\infty}\lim_{h\to\infty} S(n,h)=\lim_{h\to\infty}\lim_{n\to\infty} S(n,h)=0.$$

Therefore, from eq. (\ref{eq:zero1}), we have
\begin{eqnarray*}
\sum_{k=0}^\infty a_{2^k}=\sum_{k=1}^\infty a_k+\sum_{k=1}^\infty f(k)a_k=0
\end{eqnarray*}
and
\begin{eqnarray*}
\sum_{k=0}^\infty b_{2^k}=\sum_{k=1}^\infty b_k+\sum_{k=1}^\infty f(k)b_k=0.
\end{eqnarray*}

\indent
Since $a_1=1$, $b_1=0$ and $2^k$ is an even number for all $k$, we have
$$1-\sum_{k=1}^\infty \frac{1}{2^{kx}}\cos(ky\ln 2)=\sum_{k=0}^\infty a_{2^k}=0$$
and
$$\sum_{k=1}^\infty \frac{1}{2^{kx}}\sin(ky\ln 2)=-\sum_{k=0}^\infty b_{2^k}=0.$$
This contradicts Lemma \ref{lemf}.
Thus, the condition described above is sufficient to establish the Riemann hypothesis.
\end{proof}


\end{document}